\documentclass[12pt,twosided,reqno]{amsart}
\usepackage{amsmath}
\usepackage{amsthm}
\usepackage{amsfonts}
\usepackage{amssymb}

\theoremstyle{theorem}

\theoremstyle{definition}

\begin{document}

\vskip .5cm
\title{Center of gravity and a  characterization of parabolas}
 \vskip0.5cm

\vskip 0.5cm


\vskip 0.5cm
\author{Dong-Soo Kim$^1$, Young Ho Kim$^2$ and Sookhee Park$^3$}

\address{\newline $^{1,3}$Department of Mathematics, Chonnam National University, Gwangju
500-757,  South Korea
\\ \newline
$^{2}$Department of Mathematics, Kyungpook
National University, Daegu 702-701,  South Korea}
 \email{
dosokim@chonnam.ac.kr, yhkim@knu.ac.kr and suki112@gmail.com}

\thanks{
\newline\indent The first author was  supported by Basic Science Research Program through the
    National Research Foundation of Korea (NRF) funded by the Ministry of Education, Science and Technology (2010-0022926).
 \newline\indent The second author was   supported by Basic Science Research Program through
the National Research Foundation of Korea (NRF)
funded by the Ministry of Education, Science and Technology (2012R1A1A2042298)
 and supported by Kyungpook
National University Research Fund, 2012. }
\begin{abstract}
Archimedes determined the center of gravity of  a   parabolic section as follows.
For a parabolic section  between a parabola and any chord $AB$
on the parabola, let us  denote by $P$ the point on the parabola where the tangent is parallel to $AB$
and by $V$ the point  where  the line through $P$ parallel  to the axis of the parabola
meets the chord $AB$. Then the center $G$ of gravity  of the section
lies on $PV$ called the axis of the parabolic section with $PG=\frac{3}{5}PV$.
In this paper, we study   strictly locally convex plane curves satisfying the above
center of gravity properties. As a result, we prove that among  strictly locally convex plane curves,
those properties  characterize parabolas.
\end{abstract}

\vskip 1cm

\date{}
\maketitle

\section{Introduction}

\vskip 1cm
 Archimedes found some interesting area properties of parabolas.
Consider the region bounded by a parabola and a chord $AB$.
Let $P$ be the point on the parabola where the tangent is parallel to the chord $AB$.
The parallel line through $P$  to the axis of the parabola meets the chord $AB$ at a point
$V$.
Then,   he proved that the area of the parabolic region is $4/3$ times
the area of triangle $\bigtriangleup ABP$ whose base is the chord and the third vertex is $P$.

Furthermore, he showed that the center $G$ of gravity  of the parabolic section
lies on the segment $PV$ called the axis of the parabolic section with $PG=\frac{3}{5}PV$.
For the proofs of Archimedes, see Chapter 7 of \cite{st}.
\vskip 0.3cm
Very recently, two of the present authors showed that among strictly convex plane curves,
the above area properties of parabolic sections characterize
parabolas. More precisely, they proved as follows (\cite{kk4}).
\vskip 0.3cm

 \noindent {\bf Proposition 1.}
 Let $X$ be a strictly convex curve  in the plane
 ${\mathbb R}^{2}$. Then $X$ is a parabola if and only if  it satisfies
\vskip 0.3cm
 \noindent $(C):$
 For a point $P$ on $X$ and a chord $AB$ of $X$ parallel to the tangent of $X$ at $P$,
 the area of the region bounded by the curve and $AB$ is $4/3$ times
the area of triangle $\bigtriangleup ABP$.
 \vskip 0.3cm

Actually, in \cite{kk4}, they
 established five characterizations of parabolas, which are the converses of
well-known properties of parabolas  originally due to Archimedes
(\cite{st}).  In \cite{ks}, the first author and K.-C. Shim gave a characterization of parabolas
using area of triangles associated with a plane curve, which is a generalization of some results in \cite{kr}.
See also \cite{KKKP} for some  generalizations of results in \cite{ks}.
In \cite{kk2} and \cite{kk3}, two of the present authors
  proved  the  higher dimensional analogues of some results in \cite{kk4}.
\vskip 0.30cm
 For some characterizations of parabolas or  conic sections by  properties of tangent lines, see
\cite{kka} and \cite{kkpj}. In \cite{kk1}, using curvature function $\kappa$ and support function $h$
of a plane curve,
the first and second authors of the present paper gave a
characterization of ellipses and hyperbolas centered at the origin.

Among  the  graphs of functions,
B.  Richmond and T. Richmond  established a dozen characterizations of parabolas
using elementary techniques (\cite{r}). In their paper,  parabola means the graph of a quadratic polynomial in one
variable.
\vskip 0.30cm

In this paper, we study   strictly locally convex plane curves satisfying the above mentioned
 properties on the center of gravity.
Recall that a regular plane curve $X:I\rightarrow {\mathbb R}^{2}$ in the  plane
 ${\mathbb R}^{2}$, where $I$ is an open interval, is called {\it  convex} if, for all $s\in I$
 the trace $X(I)$ of $X$ lies entirely on one side of the closed
 half-plane determined by the tangent line at $s$ (\cite{d}).
 A regular plane curve $X:I\rightarrow {\mathbb R}^{2}$ is called {\it locally  convex} if, for each $s\in I$
 there exists an open  subinterval $J\subset I$ containing $s$ such that the curve $X|_J$ restricted to $J$
 is a convex curve.

Hereafter,  we will  say  that a locally convex curve $X$ in the  plane
 ${\mathbb R}^{2}$ is  {\it strictly locally convex} if the curve   is smooth
 (that is, of class $C^{(3)}$) and is of positive  curvature $\kappa$
 with respect to the unit normal $N$ pointing to the convex side.
 Hence, in this case we have $\kappa(s)=\left< X''(s), N(X(s))\right> >0$,
  where $X(s)$ is an arc-length parametrization of $X$.

For a smooth function $f:I\rightarrow {\mathbb R}$ defined on an open interval,
we will also say that $f$ is
{\it strictly convex} if the graph of $f$ has  positive curvature $\kappa$
with respect to the upward unit normal $N$. This condition is equivalent to the positivity of $f''(x)$ on $I$.

 \vskip 0.3cm
 First of all, in Section 2 we prove the following:

 \vskip 0.3cm

 \noindent {\bf Theorem 2.}
 Let $X$ be  a strictly locally convex plane curve  in the  plane
 ${\mathbb R}^{2}$. For a fixed point $P$ on $X$ and a sufficiently small $h>0$,
 we denote by $l$ the parallel line through $P+hN(P)$ to the tangent $t$ of the curve $X$ at $P$.
 If we let $d_P(h)$ the distance from the center $G$ of gravity  of the section of $X$ cut off by $l$
 to the tangent $t$ of the curve $X$ at $P$, then we have
   \begin{equation}\tag{1.1}
   \begin{aligned}
  \lim _{h\rightarrow 0}\frac{d_P(h)}{h}=\frac{3}{5}.
       \end{aligned}
   \end{equation}
 \vskip 0.3cm

 Without the help of Proposition 1, in Section 3
 we prove the following characterization theorem for parabolas with axis parallel to the $y$-axis,
 that is, the graph of a quadratic function.
 \vskip 0.3cm

 \noindent {\bf Theorem 3.}
 Let $X$ be  the graph of a  strictly convex function $g:I\rightarrow {\mathbb R}$  in the  $uv$-plane
 ${\mathbb R}^{2}$ with the upward unit normal $N$.
  For a fixed point $P=(u,g(u))$ on $X$ and a sufficiently small $h>0$,
 we denote by $l$ (resp., $V$) the parallel line through $P+hN(P)$ to the tangent $t$ of the curve $X$ at $P$
 (resp., the point where the parallel line through $P$ to the $v$-axis meets $l$).
  Then $X$ is an open part of a parabola with axis parallel to the $v$-axis if and only if
 it satisfies
 \vskip 0.3cm
 \noindent $(D):$ For a fixed point $P$ on $X$ and a sufficiently small $h>0$,
 the center $G$ of gravity  of the section of $X$ cut off by $l$ lies on the segment $PV$ with
   \begin{equation}\tag{1.2}
   \begin{aligned}
   PG=\frac{3}{5}PV,
      \end{aligned}
   \end{equation}
   where we denote by  $PV$ both of the segment and its length.
 \vskip 0.3cm
Note that  if $X$ is an open part of a parabola with axis which is not parallel to the $v$-axis
(for example, the graph of  $g$ given in (3.23) with $\alpha\ne 0$), then it does not satisfy Condition $(D)$.

 \vskip 0.3cm
Finally using Proposition 1, in Section 4 we prove the following characterization theorem for parabolas.
 \vskip 0.3cm
 \noindent {\bf Theorem 4.}
 Let $X$ be  a strictly locally convex plane curve  in the  plane
 ${\mathbb R}^{2}$. For a fixed point $P$ on $X$ and a sufficiently small $h>0$,
 we denote by $l$ the parallel line through $P+hN(P)$ to the tangent $t$ of the curve $X$ at $P$.
 We let $d_P(h)$ the distance from the center $G$ of gravity  of the section of $X$ cut off by $l$
 to the tangent $t$ of the curve $X$ at $P$. Then $X$ is an open part of a parabola if and only if
 it satisfies
 \vskip 0.3cm
 \noindent $(E):$ For a fixed point $P$ on $X$ and a sufficiently small $h>0$, we have
   \begin{equation}\tag{1.3}
   \begin{aligned}
   d_P(h)=\frac{3}{5}h.
      \end{aligned}
   \end{equation}
 \vskip 0.3cm

  Throughout this article, all curves are of class $C^{(3)}$ and connected, unless otherwise mentioned.
  \vskip 0.50cm

 \section{Preliminaries and Theorem 2}
 \vskip 0.50cm

 Suppose that $X$ is a strictly locally convex  curve  in the plane
 ${\mathbb R}^{2}$ with the unit normal $N$ pointing to the convex side.
 For a fixed point $P \in X$, and for a sufficiently small $h>0$, consider the  parallel line $l$  through
 $P+hN(P)$  to
 the tangent $t$ of $X$ at $P$.
 Let's denote by $A$ and $B$ the points where the line $l$ intersects the curve $X$.

 We denote by $S_P(h)$ (respectively, $R_P(h)$) the area of the region bounded by the curve $X$ and chord $AB$
 (respectively, of the rectangle with a side $AB$ and another one on the tangent $t$ of $X$ at $P$ with height $h>0$).
 We also denote by $L_P(h)$ the length  of the chord $AB$.
  Then we have $R_P(h)=hL_P(h)$.

 We may adopt  a coordinate system $(x,y)$
 of  ${\mathbb R}^{2}$ in such a way that  $P$ is taken to be the origin $(0,0)$
 and the $x$-axis is the tangent line of $X$ at $P$.
 Furthermore, we may assume that $X$ is locally
 the graph of a non-negative strictly convex  function $f: {\mathbb R}\rightarrow {\mathbb R}$.
\vskip 0.3cm

For a sufficiently small $h>0$, we have

   \begin{equation}\tag{2.1}
   \begin{aligned}
   S_P(h)&=\int _{f(x)<h}\{h-f(x)\}dx,\\
   R_P(h)&= hL_P(h)=h\int _{f(x)<h}1dx.
       \end{aligned}
   \end{equation}
The  integration is taken on the interval $I_P(h)=\{x\in {\mathbb R}| f(x)<h\}$.

On the other hand,  we also have
 $$
   S_P(h)=\int _{y=0}^{h}L_P(y)dy,
     $$
which shows that
\begin{equation}\tag{2.2}
   \begin{aligned}
 S_P'(h)=L_P(h).
        \end{aligned}
   \end{equation}

 \vskip 0.3cm
First of all, we  need the following lemma (\cite{kk4}), which is useful in this article.
\vskip 0.3cm

 \noindent {\bf Lemma 5.} Suppose that   $X$  is a strictly locally convex  curve  in the plane
 ${\mathbb R}^{2}$ with  the unit normal $N$ pointing to the convex side.  Then
 we have
  \begin{equation}\tag{2.3}
   \begin{aligned}
   \lim_{h\rightarrow 0} \frac{1}{\sqrt{h}}L_P(h)= \frac{2\sqrt{2}}{\sqrt{\kappa(P)}},
    \end{aligned}
   \end{equation}
 where $\kappa(P)$ is the curvature of $X$ at $P$ with respect to  the unit normal $N$.
 \vskip 0.3cm

 From Lemma 5, we get  a  geometric meaning of  curvature $\kappa(P)$ of a  locally strictly
 convex plane curve  $X$ at a point $P\in   M$. That is, we obtain
 \begin{equation}\tag{2.4}
   \begin{aligned}
   \kappa(P)=\lim_{h\rightarrow 0} \frac{8h}{L_P(h)^2}.
  \end{aligned}
   \end{equation}

 \vskip 0.30cm
Now, we give a proof of Theorem 2.

Let us  denote by $d_P(h)$  the distance from the center $G$ of gravity  of the section of $X$ cut off by $l$
 to the tangent $t$ of the curve $X$ at $P$.  Note that  the curve $X$ is of class $C^{(3)}$.
  If we  adopt  a coordinate system $(x,y)$
 of  ${\mathbb R}^{2}$ as in the beginning of this section,  then the curve $X$ is locally
 the graph of a non-negative strictly convex  $C^{(3)}$ function $f: {\mathbb R}\rightarrow {\mathbb R}$.
  Hence, the Taylor's formula of $f(x)$ is given by
 \begin{equation}\tag{2.5}
 f(x)= ax^2 + f_3(x),
  \end{equation}
where  $a=f''(0)/2$, and $f_3(x)$ is an $O(|x|^3)$  function.
Since  $\kappa(P)=2a>0$, we see that $a$ is positive.

  From the definition of  $d_P(h)$, we have
   \begin{equation}\tag{2.6}
   \begin{aligned}
 S_P(h)d_P(h)=\phi(h),
       \end{aligned}
   \end{equation}
where we put
   \begin{equation}\tag{2.7}
   \begin{aligned}
 \phi(h)=\frac{1}{2}\int_{f(x)<h}\{h^2-f(x)^2\}dx.
       \end{aligned}
   \end{equation}
We decompose $\phi(h)=\phi_1(h)-\phi_2(h)$ as follows:
\begin{equation}\tag{2.8}
   \begin{aligned}
 \phi_1(h)=\frac{1}{2}\int_{f(x)<h}h^2dx, \quad  \phi_2(h)=\frac{1}{2}\int_{f(x)<h}f(x)^2dx.
       \end{aligned}
   \end{equation}

It follows from the definition of $L_P(h)$ that
\begin{equation}\tag{2.9}
   \begin{aligned}
 \phi_1(h)=\frac{1}{2}h^2L_P(h).
       \end{aligned}
   \end{equation}
 Hence, Lemma 5 shows that
  \begin{equation}\tag{2.10}
   \begin{aligned}
 \lim_{h\rightarrow0}\frac{\phi_1(h)}{h^2\sqrt{h}}=\frac{\sqrt{2}}{\sqrt{\kappa(P)}}.
       \end{aligned}
   \end{equation}
 \vskip 0.3cm


 \noindent {\bf Lemma 6.} For the limit of $\phi_2(h)/(h^2\sqrt{h})$ as $h$ tends to $0$, we get
 \begin{equation}\tag{2.11}
   \begin{aligned}
 \lim_{h\rightarrow0}\frac{\phi_2(h)}{h^2\sqrt{h}}=\frac{\sqrt{2}}{5\sqrt{\kappa(P)}}.
       \end{aligned}
   \end{equation}
 \noindent {\bf Proof.}
 If we put $g(x)=f(x)^2$, then we have from (2.5)
\begin{equation}\tag{2.12}
 g(x)= a^2x^4 + f_5(x),
  \end{equation}
where  $f_5(x)$ is an $O(|x|^5)$  function. We let $x=\sqrt{h}\xi$.
Then, together with (2.5), (2.8) gives
\begin{equation}\tag{2.13}
   \begin{aligned}
 \frac{\phi_2(h)}{h^2\sqrt{h}}&=\frac{1}{2h^2\sqrt{h}}\int_{f(x)<h}g(x)dx\\
 &=\frac{1}{2h^2}\int_{a\xi^2+g_3(\sqrt{h}\xi)<1}g(\sqrt{h}\xi)d\xi,
       \end{aligned}
   \end{equation}
where we denote $g_3(\sqrt{h}\xi)=\frac{f_3(\sqrt{h}\xi)}{h}$.

Since $f_3(x)$ is an $O(|x|^3)$ function, we have for some constant $C_1$
 \begin{equation}\tag{2.14}
   \begin{aligned}
   |g_3(\sqrt{h}\xi)|\le C_1\sqrt{h}|\xi|^3.
    \end{aligned}
   \end{equation}
We also obtain from (2.12) that
 \begin{equation}\tag{2.15}
   \begin{aligned}
  \frac{ |g(\sqrt{h}\xi)-a^2h^2\xi^4|}{h^2}\le C_2\sqrt{h}|\xi|^5
    \end{aligned}
   \end{equation}
where $C_2$ is a constant.

If we let  $h \rightarrow 0$,
it follows from (2.13)-(2.15) that
\begin{equation} \tag{2.16}
  \begin{aligned}
\lim_{h\rightarrow 0} \frac{\phi_2(h)}{h^2\sqrt{h}}&=\frac{1}{2}\int_{a\xi^2<1}a^2\xi^4d\xi\\&= \frac{1}{5\sqrt{a}}.
 \end{aligned}
  \end{equation}
This completes the proof of Lemma 6.
$\square$

 \vskip 0.3cm
Together with (2.10), Lemma 6 shows  that
\begin{equation}\tag{2.17}
   \begin{aligned}
 \lim_{h\rightarrow0}\frac{\phi(h)}{h^2\sqrt{h}}=\frac{4\sqrt{2}}{5\sqrt{\kappa(P)}}.
       \end{aligned}
   \end{equation}
Since $S_P'(h)=L_P(h)$, it follows from Lemma 5 that
\begin{equation}\tag{2.18}
   \begin{aligned}
   \lim_{h\rightarrow 0} \frac{1}{h\sqrt{h}}S_P(h)= \frac{4\sqrt{2}}{3\sqrt{\kappa(P)}}.
    \end{aligned}
   \end{equation}
   Thus, together with (2.17) and (2.18), (2.6) completes the proof of Theorem 2.
 \vskip 0.50cm

\section{Proof of  Theorem 3}
 \vskip 0.5cm
In this section, we give a proof of  Theorem 3.

Let $X$ be  the graph of a  strictly convex function
$g:I\rightarrow {\mathbb R}$  in the $uv$-plane
 ${\mathbb R}^{2}$ with the upward unit normal $N$.

 For a fixed point $P=(b,c)\in X$ with $c=g(b)$,
 we denote by $\theta$ the angle between
 the normal $N(P)$ and the positive $v$-axis.
 Then we have $g'(b)=\tan \theta$ and
 $V=(b, c+wh)$ for sufficiently small $h>0$,
 where $w=\sqrt{1+g'(b)^2}=\sec \theta$.

By a  change of coordinates in the plane ${\mathbb R}^{2}$  given by
\begin{equation}\tag{3.1}
 \begin{aligned}
 u&=x\cos \theta -y\sin \theta +b,\\
  v&=x\sin  \theta +y \cos \theta +c,
   \end{aligned}
  \end{equation}
the graph $X: v=g(u), u\in I$  is represented by $X: y=f(x), x\in J$,
 $P$ by the origin and $V$ by the point $(\alpha h, h)$, where $\alpha =\tan\theta$.

Since $f(0)=f'(0)=0$, the Taylor's formula of $f(x)$ is given by
 \begin{equation}\tag{3.2}
 f(x)= ax^2 + f_3(x),
  \end{equation}
where  $a=f''(0)/2$, and $f_3(x)$ is an $O(|x|^3)$  function.
Since  $\kappa(P)=2a>0$, we see that $a$ is positive.

For a sufficiently small $h>0$, it follows from  the definition
of the center $G=(\bar{x}_P(h), \bar{y}_P(h))$ of gravity
 of the section of $X$ cut off by the parallel line $l$ through $P+hN(P)$ to the
 tangent $t$ of $X$ at $P$ that
   \begin{equation}\tag{3.3}
   \begin{aligned}
\bar{y}_P(h)S_P(h)&=\phi(h),\\
\bar{x}_P(h)S_P(h)&=\psi(h),
       \end{aligned}
   \end{equation}
where we put
   \begin{equation}\tag{3.4}
   \begin{aligned}
 \phi(h)=\frac{1}{2}\int_{f(x)<h}\{h^2-f(x)^2\}dx
       \end{aligned}
   \end{equation}
and
 \begin{equation}\tag{3.5}
   \begin{aligned}
 \psi(h)=\int_{f(x)<h}\{x(h-f(x))\}dx.
       \end{aligned}
   \end{equation}

   First of all, we get
\vskip 0.30cm
 \noindent {\bf Lemma 7.} If we let $I_P(h)=\{x|f(x)<h\}=(x_1(h), x_2(h))$, then we have
 \begin{equation}\tag{3.6}
   \begin{aligned}
 \phi'(h)&=h\{x_2(h)-x_1(h)\}=hL_P(h)\\
        \end{aligned}
   \end{equation}
   and
    \begin{equation}\tag{3.7}
   \begin{aligned}
  \psi'(h)&=\frac{1}{2}\{x_2(h)^2-x_1(h)^2\}.
       \end{aligned}
   \end{equation}
 \vskip 0.30cm
 \noindent {\bf Proof.}
      If we put $\bar{f}(x)=f(x)^2$ and $k=h^2$, then we have
 \begin{equation}\tag{3.8}
   \begin{aligned}
 2\phi(h)&=\int_{f(x)^2<h^2}\{h^2-f(x)^2\}dx\\
 &=\int_{\bar{f}(x)<k}\{k-\bar{f}(x)\}dx.
       \end{aligned}
   \end{equation}
We denote by $\bar{S}_P(k)$ the area of the region bounded by the graph of $y=\bar{f}(x)$ and
the line $y=k$. Then (3.8) shows that
\begin{equation}\tag{3.9}
   \begin{aligned}
 2\phi(h)=\bar{S}_P(k).
       \end{aligned}
   \end{equation}
It follows from (2.2) that
\begin{equation}\tag{3.10}
   \begin{aligned}
\frac{d}{dk}\bar{S}_P(k)=\bar{L}_P(k),
       \end{aligned}
   \end{equation}
where $\bar{L}_P(k)$ denotes the length of the interval $\bar{I}_P(k)=\{x\in {\mathbb R}| \bar{f}(x)<k\}$.

Since $k=h^2$,  $\bar{I}_P(k)$ coincides with the interval $I_P(h)=\{x\in {\mathbb R}| f(x)<h\}$.
Hence we get $\bar{L}_P(k)=L_P(h)$. This, together with (3.9) and (3.10) shows that
\begin{equation}\tag{3.11}
   \begin{aligned}
 \phi'(h)=hL_P(h),
       \end{aligned}
   \end{equation}
which completes the proof of (3.6).

 Finally, note that $\psi(h)$ is also given by
 \begin{equation}\tag{3.12}
   \begin{aligned}
 \psi(h)=\frac{1}{2}\int_{y=0}^h\{x_2(y)^2-x_1(y)^2\}dy,
       \end{aligned}
   \end{equation}
   which shows that (3.7) holds.

  This completes the proof of Lemma 7. $\square$
 \vskip 0.30cm
  We, now,  suppose that  $X$ satisfies Condition $(D)$. Then for each
  sufficiently small  $h>0$,
   $V=(\alpha h,h)$.
  Hence, we obtain $G=\frac{3}{5}(\alpha h,h)$. Therefore we get  from (3.3) that
     \begin{equation}\tag{3.13}
   \begin{aligned}
\frac{3}{5}hS_P(h)&=\phi(h)\\
       \end{aligned}
   \end{equation}
   and
    \begin{equation}\tag{3.14}
   \begin{aligned}
\frac{3}{5}\alpha hS_P(h)&=\psi(h).
       \end{aligned}
   \end{equation}
   It also follows from (3.13) and (3.14) that
   \begin{equation}\tag{3.15}
   \begin{aligned}
\psi(h)=\alpha \phi(h).
       \end{aligned}
   \end{equation}

   By differentiating (3.13) with respect to $h$, (3.6) shows that
  \begin{equation}\tag{3.16}
   \begin{aligned}
S_P(h)=\frac{2}{3}hL_P(h).
       \end{aligned}
        \end{equation}
   Differentiating (3.16) with respect to $h$ yields
  \begin{equation}\tag{3.17}
   \begin{aligned}
2hL_P'(h)=L_P(h).
       \end{aligned}
        \end{equation}
   Integrating (3.17) shows that
  \begin{equation}\tag{3.18}
   \begin{aligned}
L_P(h)=C(P)\sqrt{h},
       \end{aligned}
        \end{equation}
  where    $C(P)$ is a constant. Hence, it follows from Lemma 5 that
    \begin{equation}\tag{3.19}
   \begin{aligned}
L_P(h)=\frac{2}{\sqrt{a}}\sqrt{h},
       \end{aligned}
        \end{equation}
  from which we get
  \begin{equation}\tag{3.20}
   \begin{aligned}
x_2(h)-x_1(h)=\frac{2}{\sqrt{a}}\sqrt{h}.
       \end{aligned}
        \end{equation}

  Now, differentiating (3.15) and applying Lemma 7 show that
  \begin{equation}\tag{3.21}
   \begin{aligned}
x_2(h)+x_1(h)=2\alpha h.
       \end{aligned}
        \end{equation}
Hence, we get from (3.20) and (3.21) that
 \begin{equation}\tag{3.22}
   \begin{aligned}
x_1(h)=\alpha h - \frac{1}{\sqrt{a}}\sqrt{h}, \quad x_2(h)=\alpha h + \frac{1}{\sqrt{a}}\sqrt{h}.
       \end{aligned}
        \end{equation}
        Since $I_P(h)=(x_1(h), x_2(h))$, we obtain from (3.22) that the graph $X:y=f(x)$ is
        given by
           \begin{equation}\tag{3.23}
 \begin{aligned}
 f(x)= \begin{cases}
 \frac{1}{2a\alpha^2}\{2a\alpha x +1-\sqrt{4a\alpha x +1}\}, & \text{if $\alpha\ne 0,$} \\
 ax^2, & \text{if $\alpha= 0.$}
  \end{cases}
  \end{aligned}
  \end{equation}

        It follows from (3.23) that
        $X$ is an open part of the parabola defined by
     \begin{equation}\tag{3.24}
   \begin{aligned}
ax^2-2a\alpha xy +a\alpha^2y^2-y=0.
       \end{aligned}
        \end{equation}
     Note that if $\alpha\ne 0,$ the function $f(x)$ in (3.23) is defined on an interval $J$ such that
        $J\subset  ( -\infty, -1/(4a\alpha))$ or $J\subset  ( -1/(4a\alpha), \infty)$
        according  to the sign of $\alpha$.

   Finally, we use the following coordinate change from (3.1):
    \begin{equation}\tag{3.25}
     \begin{aligned}
 x&=u\cos \theta +v\sin \theta -b\cos \theta -c\sin \theta,\\
  y&=-u\sin  \theta +v \cos \theta +b\sin \theta -c\cos \theta.
  \end{aligned}
  \end{equation}
  Then, after a long calculation we see that the curve $X:v=g(u)$ is an
  open part of the parabola  determined by the following quadratic polynomial
  \begin{equation}\tag{3.26}
 \begin{aligned}
 g(u)= \begin{cases}
 aw^3(u-b)^2+\alpha(u-b) +c, & \text{if $\alpha\ne 0,$} \\
 a(u-b)^2+c, & \text{if $\alpha= 0.$}
  \end{cases}
  \end{aligned}
  \end{equation}
Note that $g(b)=c, g'(b)=\alpha$ and $\kappa(P)=2a$.
  This completes the proof of  the  if part of Theorem 3.
   \vskip 0.30cm

 It is elementary to show the only if part of Theorem 3,  or see Chapter 7 of \cite{st},
 which is originally due to Archimedes.
This completes the proof of Theorem 3.

\vskip 0.50cm
\section{Proof of  Theorem 4}
 \vskip 0.5cm
In this section, using Proposition 1, we give the proof of  Theorem 4.

 Let $X$ be  a strictly locally convex plane curve  in the  plane
 ${\mathbb R}^{2}$ with the unit normal $N$ pointing to the convex side.
 For a fixed point $P$ on $X$ and a sufficiently small $h>0$,
 we denote by $l$ the parallel line through $P+hN(P)$ to the tangent $t$ of the curve $X$ at $P$.
 We let $d_P(h)$ the distance from the center $G$ of gravity  of the section of $X$ cut off by $l$
 to the tangent $t$ of the curve $X$ at $P$.

 First, suppose that  $X$ satisfies Condition $(E)$.

 For a fixed point $P\in X$,
 we  adopt  a coordinate system $(x,y)$
 of  ${\mathbb R}^{2}$ as in the beginning of Section 2.  Then the curve $X$ is locally
 the graph of a non-negative strictly convex  $C^{(3)}$ function $f: {\mathbb R}\rightarrow {\mathbb R}$.
  Hence, the Taylor's formula of $f(x)$ is given by
 \begin{equation}\tag{4.1}
 f(x)= ax^2 + f_3(x),
  \end{equation}
where  $a=f''(0)/2$, and $f_3(x)$ is an $O(|x|^3)$  function.
Since  $\kappa(P)=2a>0$, we see that $a$ is positive.

 It follows from $(E)$ and the definition of  $d_P(h)$ that
   \begin{equation}\tag{4.2}
   \begin{aligned}
 \frac{3}{5}hS_P(h)=\phi(h),
       \end{aligned}
   \end{equation}
where we put
   \begin{equation}\tag{4.3}
   \begin{aligned}
 \phi(h)=\frac{1}{2}\int_{f(x)<h}\{h^2-f(x)^2\}dx.
       \end{aligned}
   \end{equation}
If we differentiate  $ \phi(h)$ with respect to $h$, then Lemma 7 shows that
\begin{equation}\tag{4.4}
   \begin{aligned}
 \phi'(h)=hL_P(h).
       \end{aligned}
   \end{equation}

 By differentiating both sides of (4.2) with respect to $h$, we get from (4.4) and (2.2) that
 \begin{equation}\tag{4.5}
   \begin{aligned}
 S_P(h)=\frac{2}{3}hL_P(h),
       \end{aligned}
   \end{equation}
which shows  that the curve $X$ satisfies Condition $(C)$ in Proposition 1.
Note that the argument in the proof of Proposition 1 given by \cite{kk4}
can be applied even if the curve $X$ is a strictly locally  convex plane curve.
  This completes the proof of the if part of Theorem 4.
   \vskip 0.30cm

 For  a  proof of the only if part of Theorem 4,  see Chapter 7 of \cite{st},
 which is originally due to Archimedes.
This completes the proof of Theorem 4.

 \vskip 0.50cm

\noindent {\bf Acknowledgments.}
The authors  appreciate Professor Hong-Jong Kim for suggesting to study whether among the strictly convex
plane curves, the
center of gravity properties of parabolic sections characterizes parabolas.

\vskip 1.0 cm






\vskip 0.50cm

\begin{thebibliography}{5.4}


\bibitem {d}
do Carmo, M. P., {\it Differential Geometry of
Curves and Surfaces}, Prentice-Hall, Englewood Cliffs,
NJ, 1976.

\bibitem {kka}
  Kim, D.-S.  and  Kang, S. H., {\it A characterization of conic sections},
 Honam Math. J. 33 (2011), no. 3, 335-340.

 \bibitem {KKKP}
Kim, D.-S., Kim, W., Kim, Y. H. and  Park, D. H., {\it
Area of  triangles associated with a curve II}, Bull. Korean Math. Soc., 52 (2015), No. 1, 275-286.

\bibitem {kk1}
Kim, D.-S.  and Kim, Y. H., {\it A characterization of ellipses}, Amer. Math. Monthly 114 (2007), no. 1, 66-70.

\bibitem {kk2}
Kim, D.-S.  and  Kim, Y. H., {\it Some characterizations of spheres and elliptic paraboloids}, Linear Algebra Appl.
437 (2012), 113-120.



\bibitem {kk3}
Kim, D.-S.  and  Kim, Y. H., {\it Some characterizations of spheres and elliptic paraboloids II},
 Linear Algebra Appl., 438 (2013),  1356-1364.

 \bibitem {kk4}
Kim, D.-S.  and  Kim, Y. H., {\it On the Archimedean characterization of parabolas},  Bull. Korean Math. Soc.,
50 (2013),  2103-2114.

 \bibitem {kkpj}
 Kim, D.-S.,   Park, J. H. and  Kim, Y. H., {\it Some characterizations of parabolas},
  Kyungpook Math. J. 53 (2013), no. 1, 99-104.

 \bibitem {ks}
Kim, D.-S.  and  Shim, K.-C.,  {\it Area of  triangles associated with a curve},  Bull. Korean Math. Soc.,
51 (2014), no. 3, 901-909.
\bibitem {kr}
Krawczyk, J., {\it On areas associated with a curve},
Zesz. Nauk. Uniw. Opol. Mat. 29 (1995), 97-101.

\bibitem {r}
Richmond, B. and Richmond, T.,  {\it How to recognize a parabola}, Amer. Math. Monthly 116(2009), 910-922.

\bibitem {st}
Stein, S., {\it Archimedes. What did he do besides cry Eureka?},
Mathematical Association of America, Washington, DC, 1999.

\end{thebibliography}
\end{document}